\documentclass[12pt]{article}
\usepackage{graphicx}
\usepackage{amssymb}
\usepackage{amsmath}
\numberwithin{equation}{section}

\begin{document}
\author{Lev Sakhnovich}
\date{March 29, 2005}
\textbf{Meromorphic Solutions of  Linear Differential Systems,
Painleve type functions}
\begin{center}\textbf{Lev Sakhnovich}\end{center}
\emph{735 Crawford ave., Brooklyn, 11223, New York, USA}\\ E-mail
address: lev.sakhnovich@verizon.net

\begin{center}\textbf{Abstract}\end{center}
We consider the $n{\times}n$ matrix linear differential systems in
the complex plane. We find necessary and sufficient conditions
under which these systems have meromorphic  fundamental solutions.
Using the operator identity method we construct a set of systems
which have meromorphic solutions. We prove that the well known
operator with the sine kernel generates a class of meromorphic
Painleve type functions. The fifth Painleve function belongs to
this class. Hence we obtain a new and simple proof that the fifth
Painleve function  is meromorphic.\\ \textbf{Mathematic Subject
Classification (2000)}. Primary 34M05, Secondary 34M55, 47B38.\\
\textbf{Keywords.} Strong regularity, global solution, meromorphic
solution, differential system with parameter, operator identity
method, Painleve type functions.
\section{Introduction}
\setcounter{equation}{0} Let us consider the $n{\times}n$ matrix
system of the form
\begin{equation}
\frac{dW}{dx}=A(x)W,\end{equation} where $A(x)$ is the
$n{\times}n$ matrix function. Further we assume that the
 matrix function $A(x)$ is holomorphic and single valued in a
punctured neighborhood of a point $x_{0}$.\\ Every fundamental
solution $W(x)$ of system (1.1) has the form (see
[1],[20])\begin{equation} W(x)=S(x)(x-x_{0})^{\Phi},\end{equation}
where the matrix $S(x)$ is holomorphic and single valued in  the
domain $0<|x-x_{0}|<\rho$ and $\Phi$ is a constant matrix.\\
 \textbf{Definition 1.1.} (see [1],[20]) The point $x_{0}$ is called
a \emph{ regular point of system }(1.1) if the corresponding
matrix $S(x)$ is either holomorphic in a neighborhood of $x_{0}$
or has a pole in $x_{0}$.\\
 \textbf{Definition 1.2.} We shall name the regular point $x_{0}$
of system (1.1) \emph{a strong regular} if $\Phi=0$ in formula
(1.2).\\ We use the following transformation \begin{equation}
W(x)=F(x)Y(x),
\end{equation} where \begin{equation}
F(x)=\sum_{k=\ell}^{m}f_{k}(x-x_{0})^{k},\end{equation} $f_{k}$
are constant $n{\times}n$ matrices, $\mathrm{det}F(x){\ne}0,\quad
x{\ne}x_{0}.$ Then system (1.1) takes the form \begin{equation}
\frac{dY}{dx}=B(x)Y, \end{equation}where \begin{equation}
B(x)=F^{-1}(x)A(x)F(x)-F^{-1}(x)\frac{dF}{dx}.\end{equation} The
following  theorem gives the condition of regularity .\\
\textbf{Theorem 1.1.} (Horn's theorem (see [1])\emph{The point
$x_{0}$ is regular for system $(1.1)$ if and only if there exists
 transformation $(1.3)$ such that the corresponding matrix
$B(x)$ has the form
\begin{equation} B(x)=\frac{B_{1}(x)}{x-x_{0}},\end{equation}
where $B_{1}(x)$ is holomorphic in the domain}
$0{\leq}|x-x_{0}|<\rho.$\\The conditions of  Horn's theorem are
necessary conditions of the strong regularity.
  In the present  paper we give
necessary and sufficient conditions of strong regularity.\\
Separately we consider the case when the entries of $A(x)$ are
meromorphic functions .\\ \textbf{Definition 1.3.} We name the
fundamental solution $W(x)$ of system (1.1) \emph{global strong
regular} if this solution is strong regular for all singular
points of $A(x)$.\\ It is easy to see that the global strong
solution is meromorphic.
 We apply the obtained results to the canonical differential systems [15]
with the spectral parameter $\rho$:
\begin{equation} \frac{dW(x,\rho)}{dx}=[P(x)+{\rho}Q(x)]W(x,\rho).\end{equation} We
investigate in detail the special case when $n=2,\quad P(x)=0,$
and
\begin{equation}Q(x)=\left[\begin{array}{cc}
  0 & r^{-2}(x) \\
  r^{2}(x) & 0 \\
\end{array}\right].\end{equation}
Now we shall explain the connection of system (1.8), (1.9) with
the classical second order equations. The solution of this system
\begin{equation}
U(x,\rho)=\mathrm{col}[u_{1}(x,\rho),u_{2}(x,\rho)] \end{equation}
satisfies the relations
\begin{equation}\frac{du_{1}}{dx}=i{\rho}\,r^{-2}(x)u_{2}(x,\rho),\quad
\frac{du_{2}}{dx}=i{\rho}\,r^{2}(x)u_{1}(x,\rho).\end{equation}
System (1.8), (1.9) reduces to two equations of the second order .
\begin{equation}-\frac{d}{dx}r^{2}(x)\frac{du_{1}}{dx}=\rho^{2}r^{2}(x)u_{1}(x,\rho),\end{equation}
\begin{equation}-\frac{d}{dx}r^{-2}(x)\frac{du_{2}}{dx}=\rho^{2}r^{-2}(x)u_{2}(x,\rho).\end{equation}
Let us note that equations (1.12) and (1.13) are mutually dual
[5],[8],[16] and
 play an important role in a number of
theoretical and applied problems (prediction theory [11],
vibration of a thin straight rod [3],generalized string equation
[15] ).\\
 Using the operator identity method [15] we construct
classes  $r(x)$ such that the corresponding equations (1.12) and
(1.13) have meromorphic solutions in respect to $x$ for all
$\rho$. In particular we construct a class of the rational
functions $r(x)$ with this property.\\ The operator identity
method allows to construct an analytic continuation of $r^{2}(x)$
from half-axis $(0,\infty)$ onto the complex plane. We have
applied this approach to the third and the fifth Painleve
functions.  In particular we have obtained a new
 and simple proof that the fifth Painleve function is meromorphic
( see [7]).\\ \textbf{Remark 1.1.} The global Fuchsian theory (see
[1],[20]) requires that the regularity condition be met at
infinity as well. In our approach this condition can fail for
$x=\infty$. Thus our theory can be applied to the important
examples  ( see sections 8-10) in which classical Fuchsian theory
does not work.\\ \textbf{Remark 1.2.} The meromorphic solutions of
the differential systems play an important role in the spectral
theory in the space with indefinite metric [13].
\section{ Conditions of strong regularity}
\setcounter{equation}{0} Taking into account  Horn's theorem we
begin with the matrix function $A(x)$ which can be represented in
the form
\begin{equation}
A(x)=\frac{a_{-1}}{x-x_{0}}+a_{0}+a_{1}(x-x_{0})+... ,
\end{equation}where $a_{k}$ are $n{\times}n$ matrices.
We investigate the case when $x_{0}$ is either a regular point of
$W(x)$ or a pole. Hence the following relation
\begin{equation}
W(x)=\sum_{k{\geq}m}b_{k}(x-x_{0})^{k},\quad b_{m}{\ne}0
\end{equation} is true. Here $b_{k}$ are $n{\times}n$ matrices.
We note that $m$ can be negative. From formulas (1.1), (2.1) and
(2.2) we deduce that
\begin{equation}
(k+1)b_{k+1}=\sum_{j+\ell=k}a_{j}b_{\ell},
\end{equation} where $j{\geq}-1,\quad \ell{\geq}m.$  Relation (2.3)
can be rewritten in the recurrent form \begin{equation}
[(k+1)I_{n}-a_{-1}]b_{k+1}=\sum_{j+\ell=k}a_{j}b_{\ell},\quad
k{\geq}m, \end{equation}where $j{\geq}0,\quad \ell{\geq}m.$
 When $k=m-1$ we have \begin{equation}
(mI_{n}-a_{-1})b_{m}=0. \end{equation} From relation (2.5) we
deduce the following assertion.\\ \textbf{Proposition
2.1.}(necessary condition)  \emph{If the solution of system
$(1.1)$ has form $(2.2)$ then $m$ is an eigenvalue of $a_{-1}$.}\\
We denote by M the greatest integer eigenvalue of the matrix
$a_{-1}$. Using relations (2.5) we obtain the assertion.\\
\textbf{Proposition 2.2.}(sufficient condition) \emph{If the
matrix system
\begin{equation}
[(k+1)I_{n}-a_{-1}]b_{k+1}=\sum_{j+\ell=k}a_{j}b_{\ell},
\end{equation} where $m{\leq}k+1{\leq}M$ has a solution $b_{m},
b_{m+1},...,b_{M}$ and $b_{m}{\ne}0$ then system $(1.1)$ has a
solution of form $(2.2)$.}\\
 We consider the system of equations
\begin{equation}
\frac{dY}{dx}=-YA(x), \end{equation} where $Y(x)$ has the form
\begin{equation}
Y(x)=\sum_{k{\geq}p}c_{k}(x-x_{0})^{k},\quad c_{p}{\ne}0.
\end{equation}Formulas (2.1), (2.7) and (2.8) imply that
\begin{equation}
(k+1)c_{k+1}=-\sum_{j+\ell=k}c_{\ell}a_{j},
\end{equation} where $j{\geq}-1,\quad \ell{\geq}p.$ We rewrite
relation (2.9) in the form \begin{equation}
c_{k+1}[(k+1)I_{n}+a_{-1}]=-\sum_{j+\ell=k}c_{\ell}a_{j},\end{equation}
 where $j{\geq}0,\quad \ell{\geq}p.$ In the same way as
 Propositions 2.1 and 2.2 we obtain the following results.\\
 \textbf{Proposition 2.3.}  \emph{If the solution of system $(1.1)$ has form
$(2.2)$
 then $(-p)$ is an eigenvalue of $a_{-1}$.}\\
 \textbf{Proposition 2.4.}  \emph{ Let $(-P)$ be the smallest integer
 eigenvalue of the matrix $a_{-1}$. If the matrix system
 \begin{equation}c_{k+1}[(k+1)I_{n}+a_{-1}]=-\sum_{j+\ell=k}c_{\ell}a_{j},\end{equation}
 where $p{\leq}k+1{\leq}P$ has a solution
 $c_{p}\ne0,c_{p+1},...,c_{P}$, then system $(1.1)$ has a solution of
 form  $(2.2)$.}\\
 \textbf{Remark 2.1.} If $W(x)$ is a fundamental solution of system
 (1.1), then $Y(x)=W^{-1}(x)$ is a fundamental solution of  system
 (2.7).\\
 \textbf{Proposition 2.5.} \emph{If $W(x)$ and $W^{-1}(x)$ satisfy
 relations $(1.1), (2.1)$ and $(2.7), (2.8)$ respectively, then $m$ and $-p$ are
 eigenvalues of the matrix $a_{-1}$. The corresponding matrix
 $a_{-1}$ is either scalar or has at least two different integer
 eigenvalues.}\\
 \emph{Proof.} Let the matrix $a_{-1}$ not be a scalar one. Then
 it follows from relation (2.5) that \begin{equation}
 {\mathrm{det}}\,b_{m}=0. \end{equation} Let us suppose that $a_{-1}$ doesn't
 have integer eigenvalues different from $m$. In view of
  Propositions 2.1 and 2.3
 the equality $p=-m$ is true. From relations (2.2),
 (2.8) and the
 equality $W(x)W^{-1}(x)=I_{n}$ we have $b_{m}c_{m}=I_{n}$ which
 contradicts relation (2.12).  This proves the proposition.\\
\section{Integer eigenvalues }
 We consider again  differential system (1.1) ,
where $A(x)$ has  form (2.1). Let $T$ be a constant matrix such
that
\begin{equation} T^{-1}a_{-1}T=b_{-1}, \end{equation}where
$b_{-1}$ is Jordan matrix, i.e. $b_{-1}$ has the following
structure $b_{-1}={\mathrm{diag}}(J_{1},J_{2},...,J_{s}),\quad
s{\leq}n.$ Here $J_{k}=\lambda_{k}I_{k}+H_{k},\quad
1{\leq}k{\leq}s$ and
\\
$H_{k}=\left[\begin{array}{ccccc}
  0 & 1 &0 & ... & 0 \\
  0 & 0 & 1& ... & 0 \\
  ... & ... & ... &... & ...\\
  ... & ...& ... & ... & 1 \\
  0& 0 & 0& ...&0\\
\end{array}\right].$\\
We reduce system (1.1) to the form \begin{equation}
\frac{dV}{dx}=B(x)V, \end{equation} where $W(x)=TV,\quad
B(x)=T^{-1}A(x)T$. Now we describe the "shearing" transformation (
see [20],Ch.5) which lowers the eigenvalue $\lambda_{s}$ of the
matrix $b_{-1}$ by one , while leaving the others unchanged.\\
 We denote by $q$  the order of Jordan matrix $J_{s}$ and represent
$B(x)$ in the form
\begin{equation}B(x)=\frac{1}{x-x_{0}}\left[\begin{array}{cc}
  \tilde{b}_{-1} & 0\\
  0& J_{s} \\
\end{array}\right]+\tilde{B}(x),\end{equation}where $\tilde{B}(x)$
is holomorphic at $x_{0}$ and \begin{equation}
\tilde{b}_{-1}={\mathrm{diag}}(J_{1},J_{2},...,J_{s-1}).\end{equation}
The "shearing" transformation  is defined by the relation (see
[1],[20])
\begin{equation} V=S(x)U.\end{equation} Here
\begin{equation}S(x)=\left[\begin{array}{cc}
  I_{n-q} & 0\\
  0& (x-x_{0})I_{q} \\
\end{array}\right].\end{equation} Using (3.5) we deduce that
\begin{equation} \frac{dU}{dx}=C(x)U,\end{equation}where
\begin{equation}
C(x)=S^{-1}(x)B(x)S(x)-S^{-1}(x)\frac{d}{dx}S(x) .\end{equation}It
follows from (3.3) and (3.6) that
\begin{equation}c_{-1}=\left[\begin{array}{cc}
  \tilde{b}_{-1} & 0\\
  \Gamma& J_{s}-I_{q} \\
\end{array}\right].\end{equation} It is easy to see  that the
matrix $c_{-1}$ has the same eigenvalues as $b_{-1}$ except that
the eigenvalue $\lambda_{s}$ has been decreased by unity.\\
\textbf{Theorem 3.1.} \emph{If the fundamental solution of system
$(1.1)$ is strong regular then all the eigenvalues of the
corresponding matrix $a_{-1}$ are integer.}\\ \emph{Proof.} Using
a finite number of pairs of constant and "shearing"
transformations we can reduce system (1.1) to the system
\begin{equation}
\frac{d}{dx}\tilde{W}(x)=\tilde{A}(x)\tilde{W}(x), \end{equation}
where all the integer eigenvalues of $\tilde{a}_{-1}$ coincide
with the smallest integer eigenvalue of $a_{-1}$, the  non integer
eigenvalues of $a_{-1}$ and $\tilde{a}_{-1}$ coincide. If the
fundamental solution $W(x)$  of system (1.1)  is strong regular
then the fundamental solution $\tilde{W}(x)$  of system (3.10) is
strong regular as well. If $\tilde{a}_{-1}$ has non integer
eigenvalues then according to Proposition 2.5 the matrix
$\tilde{a}_{-1}$ has at least two different integer eigenvalues.
The theorem is proved.\\
\section{Examples}
 \textbf{Example 4.1.}  V.Katsnelson and D.Volok [10] investigated the
 case when the point $x_{0}$ is a simple pole of $W(x)$ and a
 holomorphic point of the inverse matrix function $W^{-1}(x)$.
 They proved that in this case \begin{equation}
 a_{-1}^{2}=-a_{-1},\quad a_{-1}a_{0}a_{-1}=-a_{0}a_{-1}.
 \end{equation}It follows from the first of the relations
 (4.1) that
 the eigenvalues of $a_{-1}$ are equal $-1$ or $0$. From
 Proposition 2.2 we deduce the assertion.\\
 \textbf{Proposition 4.1.} \emph{Let conditions $(4.1)$ be fulfilled.
 Then
 system $(1.1)$ has a strong regular solution , where $m=-1$.}\\
\emph{Proof.} In the case under consideration we have $m=-1,\quad
M=0.$ Hence system (2.6) takes the form \begin{equation}
(I_{n}+a_{-1})b_{-1}=0,\quad
-a_{-1}b_{0}=a_{0}b_{-1}.\end{equation}Comparing the first
relations of (4.1) and (4.2) we obtain the equality
\begin{equation} b_{-1}=a_{-1}c, \end{equation} where c is an
arbitrary invertible matrix. It follows from the second relation
of (4.1)   that
\begin{equation} b_{0}=a_{0}a_{-1}c \end{equation}
satisfies the second equality of (4.2). The proposition is
proved.\\ \textbf{Example 4.2.}  Let us consider the case when
$A(x)$ has a pole of the second order. We suppose that the matrix
$A(x)$ has the form
\begin{equation}A(x)=\left[\begin{array}{cc}
  a_{11}(x) & a_{12}(x)\\
  a_{21}(x)& a_{22}(x)\\
\end{array}\right],\end{equation} where
\begin{equation}
a_{11}(x)=\alpha_{0}+\alpha_{1}(x-x_{0})+...,\end{equation}\begin{equation}
a_{22}(x)=\beta_{0}+\beta_{1}(x-x_{0})+...,\end{equation}\begin{equation}
a_{12}(x)=\gamma_{-2}(x-x_{0})^{-2}+\gamma_{-1}(x-x_{0})^{-1}+...,\end{equation}
\begin{equation}
a_{21}(x)=\mu_{2}(x-x_{0})^{2}+\mu_{3}(x-x_{0})^{3}+....\end{equation}
We introduce the matrix function \begin{equation}
\tilde{A}(x)=F^{-1}(x)A(x)F(x)-F^{-1}(x)\frac{d}{dx}F(x)
,\end{equation}where \begin{equation} F(x)=\left[\begin{array}{cc}
  \frac{1}{x-x_{0}} & 0 \\
  0 & 1 \\
\end{array}\right]\end{equation}It is easy to see that the matrix
function $ V(x)=F^{-1}(x)W(x) $ satisfies the equation
\begin{equation}
\frac{d}{dx}V=\tilde{A}(x)V.  \end{equation}It is important that
the matrix function $\tilde{A}(x)$ has a pole of the first order.
Indeed it follows from formulas (4.5)-(4.10)  that
\begin{equation}
\tilde{A}(x)=\frac{\tilde{a}_{-1}}{x-x_{0}}+\tilde{a}_{0}+...
\quad ,
\end{equation} where \\
$\tilde{a}_{-1}=\left[\begin{array}{cc}
  1 & \gamma_{-2} \\
  0 & 0 \\
\end{array}\right],$ \quad$\tilde{a}_{0}=\left[\begin{array}{cc}
  \alpha_{0} & \gamma_{-1} \\
  0 & \beta_{0} \\
\end{array}\right].$\\We shall consider the case when
\begin{equation}
V(x)=\sum_{k{\geq}m}\tilde{b}_{k}(x-x_{0})^{k},\quad
\tilde{b}_{m}{\ne}0
\end{equation} is true. Here $\tilde{b}_{k}$ are $n{\times}n$
matrices.\\ \textbf{Proposition 4.2.} \emph{Let the matrix $A(x)$
have the form defined by relations $(4.6)-(4.9)$. Then system
$(1.1)$ has a strong regular solution if and only if}
\begin{equation}
\gamma_{-2}(\alpha_{0}-\beta_{0})=\gamma_{-1}.\end{equation}
\emph{Proof.} In this case we have $m=0,\quad M=1.$  From equality
$\tilde{a}_{-1}\tilde{b}_{0}=0 $ we deduce that $\tilde{b}_{0}$
has the following form
\begin{equation} \tilde{b}_{0}=\left[\begin{array}{cc}
  -s\gamma_{-2} & -t\gamma_{-2} \\
  s & t\\
\end{array}\right].\end{equation} In view of (2.4) we have
\begin{equation}
(I_{m}-\tilde{a}_{-1})\tilde{b}_{1}=\tilde{a}_{0}\tilde{b}_{0}
\end{equation} Equation (4.17) has a solution $\tilde{b}_{1}$ if and
only if relation (4.15) is fulfilled. From Proposition 2.2 we
deduce the desired assertion.\\ \textbf{Corollary 4.1.} \emph{In
addition to the conditions of Proposition 4.2 we suppose that
$\alpha_{0}=\beta_{0}$. System $(1.1)$ has a strong regular
solution
 if and only if $\gamma_{-1}=0$.}\\
\section{Differential Systems with spectral parameter}
We consider the differential system with the parameter $\rho$:
\begin{equation} \frac{dW(x,\rho)}{dx}=[P(x)+{\rho}Q(x)]W(x,\rho),
\end{equation} where the $n{\times}n$ matrix functions $P(x)$ and
$Q(x)$ can be represented in the forms \begin{equation}
P(x)=\frac{p_{-1}}{x-x_{0}}+p_{0}+... \quad , \end{equation}
\begin{equation}
Q(x)=\frac{q_{-1}}{x-x_{0}}+q_{0}+... \quad . \end{equation}
Systems  (5.1)  play an important role in the spectral theory of
the canonical differential systems with the spectral parameter
$\rho$ (see [15]). Due to  Theorem 3.1 the following assertion is
true.\\ \textbf{Proposition 5.1.}(necessary condition) \emph{If
system $(5.1)-(5.3)$ has a strong regular fundamental solution
$W(x,\rho)$ for all $\rho$ then all the eigenvalues of the matrix
$p_{-1}+{\rho}q_{-1}$ are integer and do not depend on $\rho$.}\\
\textbf{Example 5.1.} Let $n=2$ and
$p_{-1}=\left[\begin{array}{cc}
  \lambda_{1}& \Gamma_{1} \\
  0 & \lambda_{2}\\
\end{array}\right],$
$q_{-1}=\left[\begin{array}{cc}
  0& \Gamma_{2} \\
  0 & 0\\
\end{array}\right].$\\
We assume that $\lambda_{1}$ and $\lambda_{2}$ are integer
numbers. The eigenvalues of the matrix $p_{-1}+{\rho}q_{-1}$ are
equal to $\lambda_{1}$ and $\lambda_{2}$, i.e. these eigenvalues
are integer and do not depend on $\rho$.\\ \textbf{Example 5.2.}
We consider the system \begin{equation} \frac
{d}{dx}W(x,\rho)={\rho}A(x)W(x,\rho), \end{equation} where the
matrix function $A(x)$ is defined by relations (4.6)-(4.9). We
introduce the matrix
\begin{equation}
\tilde{A}(x,\rho)={\rho}F^{-1}(x)A(x)F(x)-F^{-1}(x)\frac{d}{dx}F(x)
,\end{equation}where \begin{equation} F(x)=\left[\begin{array}{cc}
  \frac{1}{x-x_{0}} & 0 \\
  0 & 1 \\
\end{array}\right]\end{equation}The matrix
function $ V(x,\rho)=F^{-1}(x)W(x,\rho) $ satisfies the equation
\begin{equation}
\frac{d}{dx}V=[P(x)+{\rho}Q(x)]V,  \end{equation} where
\begin{equation}P(x)=-F^{-1}(x)\frac{d}{dx}F(x)=\left[\begin{array}{cc}
  \frac{1}{x-x_{0}} & 0 \\
  0 & 1 \\
\end{array}\right],\end{equation}
\begin{equation}Q(x)=F^{-1}(x)A(x)F(x).\end{equation}
It follows from (5.8) and (5.9) that \\
$p_{-1}=\left[\begin{array}{cc} 1 & 0 \\
  0 & 0 \\
\end{array}\right],$
$p_{0}=\left[\begin{array}{cc} 0 & 0 \\
  0 & 0 \\
\end{array}\right],$
$q_{-1}={\rho}\left[\begin{array}{cc} 0 & \gamma_{-2} \\
  0 & 0 \\
\end{array}\right],$
$q_{0}={\rho}\left[\begin{array}{cc} \alpha_{0} & \gamma_{-1} \\
  0 & \beta_{0} \\
\end{array}\right].$\\
 Condition (4.15) takes the form \begin{equation}
{\rho}\gamma_{-2}(\alpha_{0}-\beta_{0})=\gamma_{-1}.\end{equation}
Using Propositions 4.2 and 5.1 we obtain the following
assertion.\\ \textbf{Proposition 5.2.} \emph{Let the matrix
function $A(x)$ be defined by relations $(4.6)-(4.9)$. System
$(5.1)$ has a strong regular fundamental solution for all $\rho$
if and only if}
\begin{equation}
\alpha_{0}=\beta_{0},\quad \gamma_{-1}=0.\end{equation}
\section{Global strong regular solutions}
In sections 1-5 we investigated the strong regular solutions in a
punctured neighborhood of the singular point $x_{0}$. Now we
deduce the conditions under which the solution of system (5.1) is
strong regular for all complex $x{\ne}\infty$ and all complex
$\rho{\ne}\infty$ (global strong regular solution). It is obvious
that the global strong solution is meromorphic in $x$ and entire
in$\rho$.\\
 Let us consider
the differential system
\begin{equation}
\frac{dW}{dx}={\rho}A(x)W, \end{equation} where the $2{\times}2$
matrix function $A(x)$ has the form\\
\begin{equation}A(x)=\left[\begin{array}{cc}
0 & r(x)^{-2} \\
  r(x)^{2} & 0 \\
\end{array}\right].\end{equation}
 Here $r(x)$ is a meromorphic function in the
complex plane. We denote by $x_{k},\quad
1{\leq}k{\leq}n{\leq}\infty$ and by $y_{\ell},\quad
1{\leq}\ell{\leq}m{\leq}\infty$ the different roots of $r(x)$ and
of $r^{-1}(x)$ respectively. Proposition 5.2 implies the following
result.\\ \textbf{Theorem 6.1.} \emph{Let all the roots $x_{k}$
and $y_{\ell}$ of $r(x)$ and of $r^{-1}(x)$ respectively be
simple. The fundamental solution $W(x,\rho)$ of system $(6.1),
(6.2)$ is strong regular for all $x$ and $\rho$ if and only if}
\begin{equation}
r(x_{k})=0,\quad r^{\prime}(x_{k}){\ne}0,\quad
r^{\prime\prime}(x_{k})=0, \quad (1{\leq}k{\leq}n),\end{equation}
\begin{equation}
q(y_{\ell})=0,\quad q^{\prime}(y_{\ell}){\ne}0,\quad
q^{\prime\prime}(y_{\ell})=0, \quad (1{\leq}\ell{\leq}m),\quad
q(x)=r^{-1}(x)\end{equation}.\\The following assertions can be
proved by the direct calculation.\\ \textbf{Proposition 6.1.}
\emph{ The functions
\begin{equation} I)\quad r_{1}(x)=x,\quad II)\quad
r_{2}(x)={\mathrm{tan}}x \end{equation} satisfy all the conditions
of Theorem 6.1. The corresponding $x_{k}$ and $y_{\ell}$ are
defined by the relations:}\\ $I) n=1,m=0, x_{1}=0,\quad II)
x_{k}=k\pi,\quad y_{\ell}=\ell\pi+\frac{\pi}{2}, \quad
-\infty<k,\ell<\infty.$\\ \textbf{Proposition 6.2.} \emph{If
$r(x)$ is a polynomial  and $\mathrm{deg}r(x){\geq}2$  then $r(x)$
does not satisfy conditions $(6.3)$}.\\ \textbf{Proposition 6.3.}
\emph{The functions \begin{equation}
r_{3}(x)=\frac{x-\lambda_{1}}{x-\lambda_{2}},\quad
\lambda_{1}{\ne}\lambda_{2},\end{equation}
\begin{equation}r_{4}(x)=\frac{(x-\lambda_{1})(x-\lambda_{2})}{x-\mu_{1}}
\quad \lambda_{1}{\ne}\lambda_{2},\quad
\lambda_{1,2}{\ne}\mu_{1}\end{equation}
 do not satisfy  conditions
$(6.3),(6.4)$}.\\ Theorem 6.1 and Propositions 6.1-6.3 lead to the
following problems.\\ \textbf{Problem 6.1.} To construct
meromorphic functions $r(x)$ which satisfy  conditions (6.3) and
(6.4).\\ \textbf{Problem 6.2.} To construct rational functions
$r(x)$ which satisfy  conditions (6.3) and (6.4).\\ These problems
will be investigated in the next sections.\\
 The connection of system (6.1), (6.2) with the
classical second order equations is explained in the introduction.
\section{Operator Identity}
To solve Problems 6.1 and 6.2 we use the operator identity method
(see [15]). We introduce the operators
\begin{equation} Af=i\int_{0}^{x}f(t)dt,\quad f(x){\in}L^{2}(0,a)
\end{equation} and \begin{equation}
Sf=f(x)+\int_{0}^{a}f(t)k(x-t)dt,\end{equation} where the function
$k(x),\quad (-a{\leq}x{\leq}a)$ is continuous and \begin{equation}
k(x)=k(-x)=\overline{k(x)}.\end{equation}We use the following
operator identity \begin{equation}
AS-SA^{\star}=i(\Phi_{1}\Phi_{2}^{\star}+\Phi_{2}\Phi_{1}^{\star}).\end{equation}Here
the operators $\Phi_{1}$ and $\Phi_{2}$ are defined by the
relations
\begin{equation}
\Phi_{1}g=M(x)g,\quad \Phi_{2}g=g,\end{equation} where
\begin{equation}
M(x)=\int_{0}^{x}k(u)du+\frac{1}{2},\quad
0{\leq}x{\leq}a.\end{equation} Thus the operators $\Phi_{1}$ and
$\Phi_{2}$ map the one-dimensional space of constant numbers $g$
into $L^{2}(0,a)$. Let us consider the operator \begin{equation}
S_{\xi}f=f(x)+\int_{0}^{\xi}f(t)k(x-t)dt,\quad
f(x){\in}L^{2}(0,\xi)\quad 0{\leq}\xi{\leq}a.\end{equation}Let us
formulate
 the following results (see [13]).\\
\textbf{Theorem 7.1.} \emph{We assume:\\ There are points
$0<x_{1}<x_{2}<...$ having no limit points in $[0,a]$ and such
that the operator $S_{\xi}$ is invertible on $L^{2}(0,\xi)$ for
each $\xi{\in}[0,a)/\{x_{1},x_{2},...\}$.\\Then the matrix
function
\begin{equation} B(\xi)=\Pi^{\star}S_{\xi}^{-1}P_{\xi}\Pi,\quad
\Pi=[\Phi_{1},\Phi_{2}] \end{equation} is continuous and
nondecreasing in each of the intervals $(x_{k},x_{k+1})$. The
matrix function \begin{equation}
W(\xi,\rho)=I_{2}+i{\rho}J\Pi^{\star}S_{\xi}^{-1}P_{\xi}(I-{\rho}A)^{-1}\Pi
\end{equation} is a fundamental solution for the system
\begin{equation}
W(\xi,\rho)=I_{2}+i{\rho}J\int_{0}^{\xi}[dB(t)W(t,\rho),\end{equation}where}
\begin{equation}
J=\left[\begin{array}{cc}
  0 & 1 \\
  1 & 0 \\
\end{array}\right].\end{equation}
\textbf{Theorem 7.2.} \emph{Let $B(x)$ be constructed by $(7.8)$ .
Then $B(x)$ is continuously differentiable  in the intervals
between the singularities , and in these intervals
\begin{equation}
H(\xi)=B^{\prime}(\xi)=[h_{i}^{\star}(\xi)h_{j}(\xi)]_{1}^{2},\end{equation}
where }\begin{equation}
h_{1}(\xi)=M(\xi)+\int_{0}^{\xi}\Gamma_{\xi}(\xi,t)M(t)dt.\end{equation}
\begin{equation}
h_{2}(\xi)=1+\int_{0}^{\xi}\Gamma_{\xi}(\xi,t)dt,\end{equation}We
use here the formula
\begin{equation} S_{\xi}^{-1}f=f(x)+\int_{0}^{\xi}\Gamma_{\xi}(x,t)f(t)dt.
\end{equation}We remark that $H(\xi)$ has the special form [16]
\begin{equation}
H(\xi)=\frac{1}{2}\left[\begin{array}{cc}
  Q(\xi) & 1 \\
  1 & Q^{-1}(\xi) \\
\end{array}\right].\end{equation} It follows from relations (7.10)
and (7.11) that \begin{equation}
\frac{dW(x,\rho)}{dx}=i{\rho}JH(x)W(x,\rho).\end{equation}
Introducing $U(x,\rho)=W(2x,\rho)e^{-ix\rho}$ we reduce system
(7.17) to the form
\begin{equation}\frac{dU(x,\rho)}{dx}=i{\rho}JH_{1}(x)U(x,\rho),\end{equation}
where \begin{equation}H_{1}(x)=\left[\begin{array}{cc}
  r^{2}(x) & 0 \\
  0& r^{-2}(x) \\
\end{array}\right],\end{equation} \begin{equation}
R(x)=Q(2x).\end{equation} Let us note that obtained system (7.18),
(7.19) coincides with system (6.1), (6.2).
\section{Rational $r(x)$}
Let us consider the operator $S_{\xi}$ (see(7.7)), where $k(x)$
satisfies conditions (7.3) and is a polynomial of degree $2m$. The
kernel $k(x-t)$ can be represented in the form \begin{equation}
k(x-t)=\sum_{s=0}^{2m}x^{s}p_{s}(t),\end{equation} where
$p_{s}(t)$ are the polynomials
($\mathrm{deg}\,p_{s}(t){\leq}(2m-s))$. We introduce the matrix
\begin{equation}
A_{\xi}=[\delta_{j,s}+(x^{s},p_{j}(x))_{\xi}]_{0}^{2m}
\end{equation} and the determinant \begin{equation}
\Delta_{\xi}=\mathrm{det}A_{\xi}.\end{equation}In formula (8.2) we
used the notation
$(f,g)_{\xi}=\int_{0}^{\xi}f(t)\overline{g(t)}dt.$ The solution
$g(x,\xi)$ of the equation
\begin{equation}
S_{\xi}g=f(x) \end{equation}has the form
\begin{equation}
g(x,\xi)= f(x)-\sum_{s=0}^{2m}c_{s}(\xi)x^{s},\end{equation} where
$c_{s}(\xi)=(g,p_{s})_{\xi}.$ It follows from (8.4) and Cramer's
rule that
\begin{equation}
c_{s}(\xi)=\frac{d_{s}(\xi)}{\Delta_{\xi}},\end{equation} where
the determinant $d_{s}(\xi)$ is formed by replacing the column
under number $s$ in $\Delta_{\xi}$ by the column\\
$\mathrm{col}[(f,p_{0})_{\xi},(f,p_{1})_{\xi},...,(f,p_{2m})_{\xi}].$\\
Using (8.5) and (8.6) we have \begin{equation}
g(x)=f(x)+\int_{0}^{\xi}\Gamma_{\xi}(x,t)f(t)dt,\end{equation}where
\begin{equation}
\Gamma_{\xi}(x,t)=\frac
{1}{\Delta_{\xi}}\sum_{s=0}^{2m}D_{s}(\xi,t)x^{s}.\end{equation}Here
the determinant $D_{s}(\xi,t)$ is formed by replacing the column
under number $s$ in $\Delta_{\xi}$ by the column\\
$\mathrm{col}[p_{0}(t),p_{1}(t),...,p_{2m}(t)].$\\ The expression
\begin{equation}
\Gamma_{\xi}(0,\xi)=\frac{D_{0}(\xi,\xi)}{\Delta_{\xi}}
\end{equation}
plays an important role in our theory.\\ From (8.7) and (8.8) we
deduce the following assertion.\\
 \textbf{Theorem 8.1.} \emph{If
$k(x)$ satisfies conditions $(7.3)$ and is a polynomial then the
corresponding function $Q(\xi)$ (see $(7.16)$) is rational.}\\
\emph{Proof.} Let $f(x)=1$ . In this case formula (8.5) gives
\begin{equation}
g(x,\xi)=1+\sum_{s=0}^{2m}x^{s}R_{s}(\xi), \end{equation} where
the functions $R_{s}(\xi)$ are rational. Hence the function
$g(\xi,\xi)=h_{2}(\xi)$ is rational too. The assertion of the
theorem follows directly from the equality
\begin{equation} Q^{-1}(x)=2h_{2}^{2}(x),\end{equation}
which can be obtained  from (7.12) and (7.16). \\ We denote by\\
$\xi_{1},\xi_{2},...,\xi_{n}$ the roots of the polynomial
$\Delta_{\xi}$.\\ \textbf{Theorem 8.2.} \emph{If $k(x)$ satisfies
conditions $(7.3)$ and is a polynomial then the corresponding
matrix function $W(\xi,\rho)$ defined by relation $(7.9)$ is
entire in respect to $\rho$ and meromorphic in respect to $\xi$
with the poles in the points
 $\xi_{1},\xi_{2},...,\xi_{n}$ }.\\
 \emph{Proof.} According to (7.6) the function $M(x_)$ is a
 polynomial. Hence the function $(I-A{\rho})^{-1}M(x)$ is an entire
 function of $\rho$ and $x$. Using (7.9) and (8.8) we deduce the assertion
 of the theorem.\\
 \textbf{Remark 8.2.} We consider $Q(\xi)$ and $W(\xi,\rho)$ for all
 the complex $ \xi{\ne}\xi_{k},(1{\leq}k{\leq}n)$ and for all the complex
 $\rho$.\\
 Due to analytic continuation the equality \begin{equation}
 \frac{dW}{d\xi}=i{\rho}JH(\xi)W(\xi,\rho) \end{equation} is true
for all
 complex $ \xi{\ne}\xi_{k},(1{\leq}k{\leq}n)$  and for all the complex
 $\rho$. Here the matrix function $H(\xi)$ is defined by formula
 (7.16). Theorems 7.2 and 8.2 imply the following assertion.\\
 \textbf{Corollary 8.2} \emph{The function
$r(x)=1/[\sqrt{2}h_{2}(2x)]$ is a rational function.
 The union of the sets of the roots and the poles of $r(x)$
 coincides with the set}\\
 $\frac{1}{2}\xi_{1},\frac{1}{2}\xi_{2},...,\frac{1}{2}\xi_{n}$.\\
We shall use the  relation (see [7],Ch.4.)
\begin{equation}\frac{dg(\xi,\xi)}{d\xi}=\Gamma_{\xi}(0,\xi)g(\xi,\xi),\quad g(\xi,\xi)=h_{2}(\xi).
\end{equation}From
 relation (8.13) and Corollary 8.2 we deduce the assertion.\\
 \textbf{Corollary 8.3} \emph{If all roots of $\Delta_{\xi}$ are
 simple then the corresponding function $r(x)$ satisfies the
 conditions $(6.3), (6.4)$.}\\
 \textbf{Example 8.1.} Let us consider the
 case when
 \begin{equation}
 k(x)=x^{2}.\end{equation} In this case we have $p_{0}(t)=t^{2},
 p_{1}(t)=-2t, p_{2}(t)=1.$ Hence the determinants $\Delta_{\xi}$
 and $d_{0}(\xi)$ are defined by the relations \begin{equation}
 \Delta_{\xi}=\left|\begin{array}{ccc}
   1+\xi^{3}/3 & \xi^{4}/4 & \xi^{5}/5 \\
   -\xi^{2} & 1-2\xi^{3}/3 & -\xi^{4}/4 \\
   \xi & \xi^{2} & 1+ \xi^{3}/3  \\
 \end{array}\right|, \end{equation}
\begin{equation}
d_{0}(\xi)=\left|\begin{array}{ccc}
   -\xi^{2} & \xi^{4}/4 & \xi^{5}/5 \\
   2xi & 1-2\xi^{3}/3 & -\xi^{4}/4 \\
   -1 & \xi^{2} & 1+ \xi^{3}/3  \\
 \end{array}\right| .\end{equation}It follows from (8.15) and (8.16) that
 \begin{equation}
\Delta_{\xi}=\frac{1}{1080}\xi^{9}-\frac{1}{30}\xi^{6}+1,\end{equation}
\begin{equation}
d_{0}(\xi)=-\xi^{2}[(1-\xi^{3}/6)^{2}+\xi^{3}(1-\xi^{3}/15)/2-\xi^{3}(1-\xi^{3}/24)/5].
\end{equation}The polynomial $\Delta_{\xi}$ has nine different
roots\\ $x_{k}^{3}=6, (1{\leq}k{\leq}3),\quad
x_{k}^{3}=15+9\sqrt{5}, (4{\leq}k{\leq}6),\quad
x_{k}^{3}=15-9\sqrt{5}, (7{\leq}k{\leq}9).$\\ By the direct
calculation we prove the following assertion.\\
\textbf{Proposition 8.1.} \emph{The poles of $\Gamma_{\xi}(0,\xi)$
coincide with $x_{k},(1{\leq}k{\leq}9).$ These poles are simple
and the residues in the points $x_{k},(1{\leq}k{\leq}3)$ are equal
to $1$ and in the points $x_{k},(4{\leq}k{\leq}9)$ are equal to
$-1$.}\\ From relation (8.17) and Proposition 8.1 we deduce that
\begin{equation}
h_{2}(x)=\frac{\frac{1}{6}x^{3}-1}{\frac{1}{180}x^{6}-\frac{1}{6}x^{3}-1}.\end{equation}
Hence the corresponding function $r(x)=1/[\sqrt{2}h_{2}(2x)]$ is
rational and satisfies conditions (6.3) and (6.4) (see Problem
6.2).
\section{Exponential $r(x)$}
The following example was considered in the paper [13].\\
\textbf{Example 9.1.} Let the operator $S_{\xi}$ have the form
\begin{equation}
S_{\xi}f=f(x)+\beta\int_{0}^{\xi}[e^{i\lambda(x-t)}+e^{-i\lambda(x-t)}]f(t)dt,\end{equation}
where $\beta=\overline{\beta}{\ne}0,\quad \lambda>0.$ We find
\begin{equation}
S_{\xi}^{-1}f=f(x)-K_(x)T^{-1}(\xi)\int_{0}^{\xi}K^{\star}(t)f(t)dt,\end{equation}where
$K(x)=[e^{i\lambda{x}},e^{-i\lambda{x}}]$ and
\begin{equation}
T(\xi)=\left[\begin{array}{cc}
  \xi+{\beta}^{-1} & \lambda^{-1}e^{-i\lambda{\xi}}\mathrm{sin}\lambda{\xi} \\
  \lambda^{-1}e^{i\lambda{\xi}}\mathrm{sin}\lambda{\xi} &   \xi+{\beta}^{-1}  \\
\end{array}\right].\end{equation}
By direct calculation we have $h_{1}(x)=\frac{1}{2h_{2}(x)}$ and
\begin{equation} h_{2}(x)=\frac{u(x)}{v(x)},\end{equation}
where \begin{equation} u(x)=x+\beta^{-1}
-\lambda^{-1}\mathrm{sin}\lambda\,x,\quad v(x)=x+\beta^{-1}
+\lambda^{-1}\mathrm{sin}\lambda\,x.\end{equation} It is easy to
see that all the roots   and the poles of $h_{2}(x)$ are simple.
In the same way as Corollary 8.3. we deduce the following
assertion.\\\textbf{Proposition 9.1.} \emph{ The corresponding
function \begin{equation}r(x)=\frac{2x+\beta^{-1}
+\lambda^{-1}\mathrm{sin}2\lambda\,x}{\sqrt{2}(2x+\beta^{-1}
+\lambda^{-1}\mathrm{sin}2\lambda\,x)}.\end{equation}
 is rational and
satisfies conditions (6.3) and (6.4) }
\section{Analytic continuation,  Painleve transcendents}
1. Let us consider the operator
 \begin{equation}
    (S_\xi f)(x) = f(x) + \int_0^\xi k(x,t)f(t)\;dt
  \end{equation}
  on $L^2(0,\xi)$.\\
  \textbf{Theorem 10.1.} \emph{ Let the kernel $k(x,t),\quad
  0<x,t<\infty$ have an extension to a function $k(z,w)$ which is
  analytic as function of $z$ and $w$ in a region $G$
  such that $G$ contains the set $(0,\infty)$ and $zt{\in}G$
  whenever $z{\in}G,\quad 0<t<1.$ Then the function
  \begin{equation}\sigma(\xi,f,g)=(S_{\xi}^{-1}f,g)_{\xi}
  \end{equation}
   where $f(x)$ and $g(x)$ are entire functions of $x$, has an
   extension to a function $\sigma(z,f,g)$ which is analytic in
   $G$ except at isolated points. All finite singular points of
$\sigma(z,f,g)$ are poles}\\ \emph{Proof.} For small $\xi$, the
operator $S_\xi$
 differs from the identity operator by an operator
  of norm less than one.  Therefore $S_\xi$ is invertible for $0 \le
  \xi < \varepsilon$ for some $\varepsilon > 0$.
  For each $\xi$ in $(0,\infty)$, define $U_\xi$ from $L^2(0,1)$ to
  $L^2(0,\xi)$ by
  \begin{equation*}
    (U_\xi f)(t) = \sqrt{\frac{1}{\xi}} \,
                     f\left(\frac{ x}{\xi}\right),
    \qquad 0 < x < \xi.
  \end{equation*}
  Then $U_\xi$ maps $L^2(0,1)$ isometrically onto $L^2(0,\xi)$,
  and
  \begin{equation*}
    (U_\xi^{-1} g)(x) = \sqrt{\xi}
                g(t\xi)
    \qquad 0 < t < 1.
  \end{equation*}
  Hence $\tilde{S_{\xi}}=U_\xi^{-1} S_\xi U_\xi$ is a bounded operator on
  $L^{2}(0,1)$ given by
  \begin{equation}
   \tilde{S_{\xi}}f(x)=
   f(x) + \xi\int_0^1
               k(\xi x, \xi t)
               f(t) \; dt .
  \end{equation}
  Clearly $S_\xi$ is invertible if and only if $\tilde{S_{\xi}}$ is invertible.  Write
\begin{equation}\tilde{S_{\xi}}=I+T_{\xi}.\end{equation}
  The assumptions of the theorem allow us to define an operator $T(z)$ on
  $L^2(0,1)$ by
  \begin{equation}
    T(z)f=z \int_0^1
               k(zx,zt)
               f(t) \; dt  .
\end{equation}
  The operator $T(z)$ is compact and depends holomorphically on~$z$,
  and $T(z)$ agrees with the operator $T(\xi)$ defined by
 (10.4) when $z = \xi$ is a point of $(0,\infty)$.  Since $I +
  T(\xi)$ is invertible for small positive~$\xi$,  $I+T(z)$ is
  invertible except at isolated points of $G$ (see Kato [9]
  Theorem 1.9 on p. 370) in which $[I+T(z)]^{-1}$ has the poles.
  Hence the function $([I+T(z)]^{-1}x^{m},x^{n})_{1}$
is meromorphic if $m$ and $n$ are non-negative integers. The
assertion of  the theorem follows  from the relation
\begin{equation}
(S_{\xi}^{-1}x^{m},x^{n})_{\xi}=\xi^{m+n+1}({\tilde{S}}^{-1}_{\xi}x^{m},x^{n})_{1}.
\end{equation}
\textbf{Remark 10.1.}  The arguments close to Theorem 10.1 are
contained in the article [13]\\ The following kernels satisfy the
conditions of Theorem 10.1:\\
\begin{equation}
 k_{1}(x,t)=\gamma\frac{\mathrm{sin}\pi(x-t)}{\pi(x-t)},\qquad
\gamma=\overline{\gamma},\end{equation}
\begin{equation}
k_{2}(x,t)=\gamma\frac{Ai(x)Ai^{\prime}(t)-Ai(t)Ai^{\prime}(t)}{x-t},\qquad
\gamma=\overline{\gamma},
\end{equation}where $Ai(x)$ is the Airy function.\\
 Let us introduce the functions
$\phi(x)=J_{\alpha}(\sqrt{x}),\quad \psi(x)=x\phi^{\prime},\quad
x{\geq}0$  and the kernel
 \begin{equation}
 k_{3}(x,t)=\gamma\frac{\phi(x)\psi(t)-\phi(t)\psi(x)}{x-t},\qquad
 \gamma=\overline{\gamma},
 \end{equation} where $J_{\alpha}(x)$ is the Bessel function  of
 order $\alpha\quad (\alpha>-1)$.\\
 \textbf{Remark 10.2.} The sine-kernel $k_{1}(x,t)$, the
 Airy-kernel $k_{2}(x,t)$ and the Bessel-kernel $k_{3}(x,t)$ play
 an important role in the random matrix theory (see [4], [12], [17-19].)\\
 \textbf{Remark 10.3.} The region $G$ in the cases $k_{1}(x,t)$ and
 $k_{2}(x,t)$  is the complex plane. The region $G$ in the case
 $k_{3}(x,t)$
is the complex plane cut by the half-axis  $[0,\infty)$.\\ 2.
\textbf{Example 10.1 ( fifth Painleve transcendent)}\\ Let us
consider the operator\begin{equation}
S_{t}f=f(x)+\gamma\int_{-t}^{t}k(x-u)f(u)du,\quad
f(u){\in}L^{2}(-a,a),
\end{equation} where $|t|{\leq}a, \quad \gamma=\overline{\gamma} $ and \begin{equation}
k(x)=\frac{\sin{x\pi}}{x\pi}.\end{equation}The operator $S_{t}$ is
invertible (see [4],p.167), when $|\gamma|{\leq}1$. Hence we have
\begin{equation}
S_{t}^{-1}f=f(x)+\int_{-t}^{t}\Gamma_{t}(x,u,\gamma)f(u)du,\quad
f(u){\in}L^{2}(-t,t),\end{equation}where the kernel
$\Gamma_{t}(x,u,\gamma)$ is jointly continuous to the variables
$x,t,u, \gamma$. Together with the operator $S_{t}$ we shall
consider the operator
\begin{equation}
\tilde{S}_{2t}f=f(x)+\gamma\int_{0}^{2t}k(x-v)f(v)dv,\quad
f(u){\in}L^{2}(0,2t).\end{equation}  The operator
\begin{equation}
U_{t}f(x)=f(u+t)\end{equation}  maps unitarily the space
$L^{2}(0,2t)$ onto $L^{2}(-t,t)$. It is easy to see that
\begin{equation}U_{t}^{-1}S_{t}U_{t}f=
\tilde{S}_{2t}f.\end{equation}
 In view of (10.12) and (10.15) we have
\begin{equation}
\tilde{S}_{2t}^{-1}f=f(x)+\int_{0}^{2t}\tilde{\Gamma}_{2t}(x,u,\gamma)f(u)du,\quad
f(u){\in}L^{2}(0,2t),\end{equation}where \begin{equation}
\tilde{\Gamma}_{2t}(x,y,\gamma)=\Gamma_{t}(x-t,y-t,\gamma).\end{equation}It
follows from (10.17) that
\begin{equation}\tilde{\Gamma}_{2t}(2t,2t,\gamma)=\Gamma_{t}(t,t,\gamma),\quad
\tilde{\Gamma}_{2t}(2t,0,\gamma)=\Gamma_{t}(t,-t,\gamma).\end{equation}Now
we consider the case when $\gamma=-1$. For brevity we omit the
parameter $\gamma=-1$ in the notation $\Gamma_{t}(x,u,-1)$.
 Following C.Trace
and H.Widom [17] we introduce the functions \begin{equation}
r(t)=e^{it\pi}+\int_{-t}^{t}\Gamma_{t}(t,u)e^{iu\pi}du
\end{equation} and
\begin{equation}
q(t)=e^{it\pi}+\int_{0}^{t}\tilde{\Gamma}_{t}(t,u)e^{iu\pi}du
\end{equation}
Relations (10.17) and (10.19), (10.20) imply that
\begin{equation}
q(2t)=r(t)e^{it\pi}.\end{equation}
 We use the following relation (see [17])
\begin{equation}
\frac{d}{dt}[tR(t,t)]=|r(t)|^{2},\end{equation} where
$R(t,t)=\Gamma_{t}(t,t)$. From (10.21) and (10.22) we have
\begin{equation}
tR(t,t)=\frac{1}{2}\int_{0}^{2t}|q(v)|^{2}dv.\end{equation}  To
prove the relation
\begin{equation}tR(t,t)=\frac{1}{2}(\tilde{S}_{2t}^{-1}e^{u\pi},e^{u\pi})_{2t}
\end{equation} we use the notion of the triangular factorization
(see[7], Ch.4; [14]).\\
 \textbf{Definition 10.1.} The positive operator $S$
acting in $L^{2}(0,a)$ admits the \emph{triangular factorization}
if it can be represented in the form
\begin{equation}S=S_{-}S_{-}^{\star}.\end{equation}Here
$Q_{\xi}S_{-}^{\pm1}=Q_{\xi}S_{-}^{\pm1}Q_{\xi} ,$ where
  $Q_{\xi}=I-P_{\xi},\quad P_{\xi}f=f(x) ,\quad 0{\leq}x<\xi$ and
 $P_{\xi}f=0 ,\quad \xi{\leq}x{\leq}a,\quad f(x)\in L_{k}^{2}(0,a)$.\\
Using M.G.Krein result (see [7], Ch.4) on the triangular
factorization of the operator  $S$ with continuous kernel we
obtain the assertion.\\ \textbf{Theorem 10.2.} \emph{The operator}
\begin{equation}
Sf=f(x)-\frac{1}{\pi}\int_{0}^{a}\frac{{\mathrm{sin}}\pi(x-t)}{x-t}f(t)dt
\end{equation} \emph{admits  triangular factorization $(10.25)$ and}
\begin{equation}
S_{-}^{-1}f=f(v)+\int_{0}^{v}\tilde{\Gamma}_{v}(v,u)f(u)du.
\end{equation}
Hence formula  (10.20) can be written in the form \begin{equation}
q(x)=S_{-}^{-1}e^{iu\pi}.\end{equation}
 \textbf{Remark 10.4.}
Representation (10.28) of $q(x)$ which contains the factorizing
operator $S_{-}$ plays an essential role in our approach.\\ We use
the notations \begin{equation}
D(\xi)=\mathrm{det}\,\tilde{S}_{\xi},\end{equation}
\begin{equation}
\sigma(x)=\frac{x}{\pi}D^{\prime}(\frac{x}{\pi})/D(\frac{x}{\pi}).\end{equation}
It is known (see[17]) that \begin{equation}
\sigma(x)=-2tR(t,t),\quad x=2{\pi}t.\end{equation}  Relations
(10.24) and (10.31) imply that
\begin{equation}\sigma(x)=-(\tilde{S}_{2t}^{-1}e^{iu\pi},e^{iu\pi})_{2t},\quad
x= 2{\pi}t.\end{equation} We note that the function $\sigma(x)$ is
the fifth Painleve transcendent (see [17]).Using Proposition 10.1
and relation (10.31) we have obtained  the new proof of the
following well-known fact (see[6]). \\ \textbf{Corollary 10.1.}
\emph{The fifth Painleve transcendent $\sigma(\xi)$ can be
extended to the meromorphic function $\sigma(z)$.}\\ The function
$\sigma(\xi)$ is a solution of the Painleve equation ($P_{5}$ in
the sigma form, see [17])
\begin{equation}
(\xi\sigma^{\prime\prime})^{2}+4(\xi\sigma^{\prime}-\sigma)(\xi\sigma^{\prime}-\sigma+{\sigma^{\prime}}^{2})
=0.
\end{equation}
\textbf{Proposition 10.1.} \emph{All the poles $z_{k}$ of
$\sigma(z)$ are simple with residues  $z_{k}.$}\\ \emph{Proof.}
Looking at the Laurent expansion of  $\sigma(z)$ at the poles
$z_{k}$ we observe by  (10.33) that  the principal term of
$\sigma(z)$has to be $z_{k}/(z-z_{k})$. The proposition is proved.
\\
 According to (10.4) and (10.5) the function $D(\xi)$ can be
extended to the entire function $D(z)$. From Proposition 10.1. and
relation (10.29) we obtain the assertions.\\ \textbf{Corollary
10.2.} \emph{All the zeroes of $D(z)$ are simple.}\\
\textbf{Corollary 10.3.} \emph{All the eigenvalues of $T(z)$ are
simple.}\\ 3. \textbf{Example 10.2. (Painleve type functions)}\\
Let us consider the operator $S_{\xi}$ of form (10.1), where
$k(x,t)=k_{1}(x,t)$. We introduce the functions \begin{equation}
\sigma_{1}(\xi,\gamma,\lambda)=(S_{\xi}^{-1}f,g)_{\xi},\end{equation}
where $f(x)=g(x)=e^{ix\lambda},\quad \lambda=\overline{\lambda}.$
Using Theorem 10.1 and Corollary 10.3 we obtain the following
assertion.\\ \textbf{Proposition 10.2.} \emph{The function
$\sigma_{1}(\xi,\gamma,\lambda)$ can be extended to the
meromorphic function $\sigma_{1}(z,\gamma,\lambda)$, all the poles
of $\sigma_{1}(z,\gamma,\lambda)$ are simple.}\\
\textbf{Definition 10.2.} We call the functions
$\sigma_{1}(z,\gamma,\lambda)$ \emph{the Painleve type
functions.}\\ We note that the  fifth Painleve transcendent
$\sigma(z)$ is connected with the functions of form (10.34) by the
relation
\begin{equation}\sigma(z)=-\sigma_{1}(z/\pi,-1,\pi).\end{equation}
4. We separately consider the function \begin{equation}
\sigma_{2}(z,\gamma)=\sigma_{1}(z,\gamma,0).
\end{equation}  It follows from (10.34) and (10.36) that
 \begin{equation}\sigma_{2}(z)=(S_{\xi}^{-1}1,1)_{\xi},\quad
\xi>0,\end{equation}where the operator $S_{\xi}$ and the kernel
$k(x,y)$ are defined by  relations (10.1) and (10.7) respectively.
We introduce the operators of form (10.10) with the kernels
$k(x,y)$ and
\begin{equation}
 k_{\pm}(x,y)=\frac{1}{2}[k(x,y){\pm}k(-x,y)].\end{equation}
 Let us denote the Fredholm determinants corresponding
 to $k(x,y),\quad k_{+}(x,y)$ and $k_{-}(x,y)$ by
 $D(\gamma,t),\quad D_{+}(\gamma,t)$ and $D_{-}(\gamma,t)$
 respectively. We use the following relations (see [12], Ch.21.)
 \begin{equation}
 D(\gamma,t)=D_{+}(\gamma,t)D_{-}(\gamma,t),\end{equation}
\begin{equation}
 \frac{D_{-}(\gamma,t)}{D_{+}(\gamma,t)}=1+\int_{-t}^{t}\Gamma_{t}(t,y,\gamma)dy.
 \end{equation} It follows from (10.17) and (10.40) that
\begin{equation}
 \frac{D_{-}(\gamma,t)}{D_{+}(\gamma,t)}=1+\int_{0}^{2t}\tilde{\Gamma}_{2t}(t,y,\gamma)dy.
 \end{equation} In view of (7.14) and (10.40)the relation
\begin{equation}
\frac{D_{-}(\gamma,t)}{D_{+}(\gamma,t)}=h_{2}(2t)\end{equation}is
true. Using formulas  (10.3)-(10.5) we deduce the assertion\\
\textbf{Proposition 10.3.} \emph{The functions $D_{-}(x,t)$ and
$D_{-}(x,t)$ can be extended to the entire functions $D_{-}(z,t)$
and $D_{-}(z,t)$ respectively}.\\Hence we have\\ \textbf{Corollary
10.4.} \emph{The function $h_{2}(2t)$ can be extended to the
meromorphic function }
\begin{equation}
 h_{2}(2z)=\frac{D_{-}(\gamma,z)}{D_{+}(\gamma,z)}.\end{equation}According
to representation (7.9) and Theorem 10.1 equations (1.12) and
(1.13) have the strong regular solutions $u_{1}(x,\rho)$ and
$u_{2}(x,\rho)$ respectively when
\begin{equation} r^{-2}(x)=2h_{2}^{2}(2x).\end{equation} From
Theorem 6.1, Corollary 10.4 and relation (10.44) we obtain \\
\textbf{Corollary 10.5.} \emph{The function
\begin{equation}
 r(z)=\frac{D_{-}(\gamma,z)}{D_{+}(\gamma,z)}\end{equation}
satisfies conditions $(6.3),(6.4)$.}\\ 4. \textbf{Example 10.3.
(third Painleve function).}\\ Let us consider the operator (10.1),
where $k(x,t)$ is the Bessel kernel.\\ \textbf{Proposition 10.4.}
\emph{The operator $S_{\xi}$ defined by relations $(10.1)$ and
$(10.9)$ is invertible on $L^{2}(0,\xi)$, when
$|\gamma|{\leq}1$.}\\ \emph{Proof.} The kernel $k(x,t)$ has the
form $k(x,t)={\gamma}K(x,t),$ where
\begin{equation}
K(x,t)=\frac{1}{4}\int_{0}^{1}\phi(xs)\phi(ts)ds, \quad
\phi(x)=J_{\alpha}(\sqrt{x}).\end{equation} The operator
$T_{\xi}=S_{\xi}-I$ has the kernel $k(x,t)$ and is self-adjoint.
It follows from (10.1) and (10.46) that
\begin{equation}
(T_{\xi}f,f)=\frac{\gamma}{4}\int_{0}^{1}|\int_{0}^{\xi}J_{\alpha}(\sqrt{xs})f(x)dx|^{2}ds.\end{equation}
The last relation can be written in the form \begin{equation}
(T_{\xi}f,f)={\gamma}\int_{0}^{1}|F(s)^{2}ds,\end{equation}where
\begin{equation}
F(s)=\int_{0}^{\sqrt{\xi}}\sqrt{xs}J_{\alpha}(sx)f(x^{2})\sqrt{2x}dx=
\end{equation}The
Hankel transformation (10.49) is unitary. So we have
\begin{equation}|(T_{\xi}f,f)|{\leq}|\gamma|\int_{0}^{\infty}|F(s)|^{2}ds=
|\gamma|\int_{0}^{\xi}|f(x)|^{2}dx.\end{equation}Hence
$||T_{\xi}||{\leq}|\gamma|$. If  $|\gamma|<1$ then the operator
$S_{\xi}$  is invertible. We shall consider separately the case
when $\gamma={\pm}1$. Let us assume that $||T_{\xi}||=1.$ In this
case we have for some $f$ the equality \begin{equation}
T_{\xi}f={\pm}f,\quad ||f||{\ne}0.\end{equation}From relations
(10.48), (10.50) and (10.51) we deduce that \begin{equation}
F(s)=0,\quad s>1.\end{equation}But the function $F(s)$ is analytic
when $Res>0$. Hence the equality $F(s)=0, (s>0)$  is true. It
means that  $||f(x)||=0$ . We have obtained a contradiction, i.e.
$||T_{\xi}||<1$. The proposition is proved.\\The operator
$S_{\xi}^{-1}$ has the form
\begin{equation}
S_{\xi}^{-1}f=f(x)+\int_{0}^{\xi}\Gamma_{\xi}(x,t)f(t)dt.\end{equation}
  We consider the functions
\begin{equation}
q(\xi)=\phi(\xi)+\int_{0}^{\xi}\Gamma_{\xi}(\xi,t)\phi(t)dt,\end{equation}
We shall use the following relations (see [18])
\begin{equation}[sR(s)]^{\prime}=\frac{1}{4}q^{2}(s)\end{equation}
\begin{equation}R(t)=-\frac{d}{dt}\mathrm{logdet}S_{t}.\end{equation}
Using M.G.Krein result (see [7],  Ch.4) on the triangular
factorization of the operator  $S$ with continuous kernel we
obtain the assertion.\\ \textbf{Theorem 10.3.} \emph{The operator
$S_{a}$ defined by $(10.1),(10.9)$ when $\alpha{\geq}0$ admits
triangular factorization $(10.25)$ and}
\begin{equation}
S_{-}^{-1}f=f(v)+\int_{0}^{v}\Gamma_{v}(v,u)f(u)du.
\end{equation}
Formula  (10.54) can be written in the form \begin{equation}
q(x)=S_{-}^{-1}\phi.\end{equation}We introduce the notation
\begin{equation}\sigma(s)=sR(s,s).\end{equation}
Relations (10.55) and (10.58) imply that
\begin{equation}\sigma(x)=\frac{1}{4}(S_{\xi}^{-1}\phi,\phi)_{\xi}.
\end{equation}Further we consider only the important case when $\gamma=-1,
\alpha{\geq}0$.\\
 We note that in this case the function $\sigma(x)$ is
the third  Painleve transcendent (see [17]). Using Proposition
10.1 and relation (10.60) we obtain  the following fact (see[7]).
\\ \textbf{Corollary 10.6.} \emph{The third Painleve transcendent
$\sigma(\xi)$ can be extended to the  function $\sigma(z)$ which
is analytic in $G$ except at isolated points. All finite singular
points in $G$ are poles. (The domain $G$ is defined in Remark
10.3.)}\\ The function $\sigma(\xi)$ is a solution of the Painleve
equation ($P_{3}$)
\begin{equation}
(\xi\sigma^{\prime\prime})^{2}+\sigma^{\prime}(\sigma -
\xi\sigma^{\prime})(4\sigma^{\prime}-1)-{\alpha}^{2}{\sigma^{\prime}}^{2}
=0.
\end{equation}
\textbf{Proposition 10.5.} \emph{All the poles $z_{k}$ of
$\sigma(z)$ are simple with residues  $z_{k}.$}\\ \emph{Proof.}
Looking at the Laurent expansion of  $\sigma(z)$ at the poles
$z_{k}$ we observe by  (10.61) that  the principal term of
$\sigma(z)$ has to be $z_{k}/(z-z_{k})$. The proposition is
proved. \\
\begin{center}\textbf{References}\end{center}
 1. Bieberbach
L.,Theorie der Gew\"{o}nliche Differentialgleichungen,
Springer-Verlag, Berlin, 1965.\\ 2.Coddington E.A.,Levinson
N.,Theory of Ordinary Differential Equations, McGraw-Hill Book
Company, New York , 1955.\\ 3.Coleman C.F., McLaughlin J.R.,
Solution of the Inverse Spectral Problem for an impedance with
integrable derivative, Comm.Pure Appl.Math.56,145-184, 1993\\
4.Deift P.A.,Its A.R., Zhou X., A Riemann-Hilbert Approach  to
Asymptotic Problems Arising in the Theory of Random Matrix Models,
and also in the Theory of Integrable Statistical Mechanics, Annals
of Math.146, 149-235, 1997.\\ 5.Dym H., Sakhnovich L.,On dual
Canonical Systems and Dual Matrix String Equations,Operator
Theory,v.123,207-228,2001. \\ 6.Gohberg I.,Krein M.G., Theory and
Applications of Volterra Operators in Hilbert Space,Amer.
Math.Soc.Providence, 1970.\\ 7.Gromak V.I., I.Laine, S. Shimomura,
Painleve Differential Equations in the Complex Plane,  de Gruyter
Studies in Mathematics, 28, Berlin, New York, 2002.\\ 8.Kac I.S.,
Krein M.G.,On the Spectral Function of the String,Amer.Math. Soc.
Translation 103,1-18,1974.\\ 9.Kato T.,Perturbation Theory for
Linear Operators,Springer- Verlag,Berlin,1976.\\ 10.Katsnelson V.
and Volok D., Rational Solutions of theSchlesinger System and
Isoprincipal Deformations of Rational Matrix Functions II.
Preprint, 2004.\\ 11.Krein M.G. On Main Approximation Problem of
Extrapolation Theory and Filtration of Stationary Stochastic
Processes, Dokl. Akad. Nauk SSSR 94,No.1,13-16,1954\\ 12.Mehta
M.L., Random Matrices,Academic press, Boston,1991.\\ 13.Rovnyak J.
and Sakhnovich L.A., Inverse Problems for Canonical Differential
Systems with Singulariries, Preprint, 2005.\\ 14.Sakhnovich L.A.,
Factorization of Operators in $L^{2}(a,b)$, Functional Anal.
Appl.13,187-192 (Russian) 1979.\\ 15. Sakhnovich L.A., Spectral
Theory of Canonical Differential Systems. Method of Operator
Identities, Operator Theory Advances and Appl., Birkh\"{a}user,
v.107, 1999.\\ 16. Sakhnovich L.A., On Reducing the Canonical
System to the Two Dual Differential Systems, J. Math. Anal.
Appl.255 ,No.2,499-509, 2001.\\ 17.Tracy C.A. and Widom
H.,Introduction to Random Matrices,Springer Lecture Notes in
Physics 424,103-130,1993.\\ 18.Tracy C.A. and Widom H.,Level
spacing distribution and the Bessel kernel, Common. Math. Phys.
161,289-309,1994.\\ 19.Tracy C.A. and Widom H.,Level spacing
distribution and the Airy kernel, Common. Math. Phys.
159,151-174,1994.\\ 20. Wasow W.,Asymptotic Expansions for
Ordinary Differential Equations, Pure and Appl. Math., v.14, 1965.

\end{document}